\documentclass{ifacconf}

\usepackage{graphicx}     
\usepackage{natbib}

\usepackage{wasysym}

\def\dashfill{\cleaders\hbox to .5em{\rule{.4ex}{.4pt}}\hfill}
\newcommand\dashline[1]{\hbox to #1{\dashfill\hfil}}

\newcommand\vdashline{\, \brokenvert \,}

\usepackage{color}

\usepackage{array}       

\newcolumntype{P}[1]{>{\centering\arraybackslash}p{#1}}
\newcolumntype{M}[1]{>{\centering\arraybackslash}m{#1}}
       
\begin{document}
\begin{frontmatter}

\title{Stacked adaptive dynamic programming with unknown system model} 


\author[TUC]{Pavel Osinenko} 
\author[TUC]{Thomas G{\"o}hrt}
\author[TUC]{Grigory Devadze} 
\author[TUC]{Stefan Streif}

\address[TUC]{Laboratory for Automatic Control and System Dynamics; \\ Technische Universit\"at Chemnitz, 09107 Chemnitz, Germany}

\begin{abstract}
Adaptive dynamic programming is a collective term for a variety of approaches to infinite-horizon optimal control. Common to all approaches is approximation of the infinite-horizon cost function based on dynamic programming philosophy. Typically, they also require knowledge of a dynamical model of the system. In the current work, application of adaptive dynamic programming to a system whose dynamical model is unknown to the controller is addressed. In order to realize the control algorithm, a model of the system dynamics is estimated with a Kalman filter. A stacked control scheme to boost the controller performance is suggested. The functioning of the new approach was verified in simulation and compared to the baseline represented by gradient descent on the running cost.
\end{abstract}


\end{frontmatter}
\section{Introduction and problem statement} \label{sec:intro} 

Consider a discrete-time system in the form:

\begin{equation} \label{eq:sys-general}
	x_{k+1} = f(x_k, u_k),
\end{equation}

where $k$ denotes the time step, $u$ is called \emph{control action}, $x \in X$ is called \emph{state} and $X$ is the \emph{state space}, $f$ is the \emph{state transition function}. It is assumed that $u$ and $x$ are vectors in Euclidean spaces. A concrete function $v(x_k)$ in place of $u_k$ in (\ref{eq:sys-general}) is called \emph{control policy}. Infinite--horizon optimal control is the goal of solving the following optimization problem:

\begin{equation} \label{eq:infinite-hor-opt}
	\min_v \sum_{i=k}^{\infty} r( x_i, v(x_i) ), \forall x_k \in X.
\end{equation}

Here, $r$ describes the so-called \emph{running cost} or \emph{reward}. The function 

\begin{equation} \label{eq:val-fnc}
	J_v(x_k) = \sum_{i=k}^{\infty} r( x_i, v(x_i) ), x_k \in X,
\end{equation}

is also called \emph{value function} or \emph{cost-to-go} for the policy $v$.

Rewriting (\ref{eq:val-fnc}) as

\begin{equation}
	J_v(x_k) = r(x_k, v(x_k)) + J_v(x_{k+1}), x_k \in X
\end{equation}

leads to the famous Bellman equation:

\begin{equation} \label{eq:Bellman-eqn}
	J^*(x_k) = \min_v \left\{ r(x_k, v(x_k) ) + J_v( x_{k+1} ) ) \right\}, \forall x_k \in X.
\end{equation}

The Bellman's optimality principle \citep{Bellman1957-DP, Blackwell1965-DP}, inspired to some extent by \citet{Wald1947-decision-making}, essentially states the following:

\begin{equation} \label{eq:DP-principle}
	J^*(x_k) = \min_v \left\{ r(x_k, v(x_k) ) + J^*( x_{k+1} ) \right\}, \forall x_k \in X.
\end{equation}

The optimal control policy is thus determined by:

\begin{equation} \label{eq:opt-ctrl-policy}
	v^*(x_k) = \arg \min_v \left\{ r(x_k, v(x_k) ) + J^*( x_{k+1} ) \right\}, \forall x_k \in X.
\end{equation}

This is the core principle of dynamic programming. Unfortunately, solving (\ref{eq:DP-principle}) and (\ref{eq:opt-ctrl-policy}) is intractable in a number of applications since it requires exploring the whole state space. Reinforcement learning \citep{Sutton1998-RL} and adaptive dynamic programming (ADP) were designed to ease this problem of dynamic programming by using approximations of the value function \citep{Lewis1995-opt-ctrl}. Giving up some precision due to approximation made the problem (\ref{eq:DP-principle}) feasible. As approximators, neural networks attracted much popularity \citep{Werbos1990-menu,Werbos1992-ADP,Bertsekas1995-NDP}. For recent surveys on ADP, refer, for example, to \citet{Abu2005-ADP}, \citet{Balakrishnan2008-ADP-survey}, \citet{Lewis2009-ADP}, \citet{Ferrari2011-ADP-survey} and \citet{Lewis2013-RL-ADP}. Convergence analyses of ADP may be found in \citep{Al-Tamimi2008-VI-ADP,Heydari2014-better-prf-ADP,Liu2014-PI-ADP-conv-prf}. 

Adaptive dynamic programming is not a single, but rather a variety of methods. The common feature of all them is that they treat (\ref{eq:DP-principle}) and (\ref{eq:opt-ctrl-policy}) in an iterative manner. In the so-called \emph{value iteration}, one starts with an arbitrary control policy and a continuous positive-definite value function. Then, one updates the value function using this policy. In the next step, the new policy is computed by optimizing the previous value function and so on. \emph{Policy iteration}, in contrast, starts with a policy that stabilizes the system and yields a finite cost-to-go. Using this policy, one finds the value function that satisfies the Bellman equation. Optimizing this value function yields the next policy and so on. In the so-called \emph{dual learning}, one essentially performs iterations with the gradients of the value function. A survey and details of all the methods may be found in \citep{Lewis2009-ADP}. Common to all of the methods, the system dynamics and the gradient of the state with respect to the control must be known. The latter is not required in the so-called \emph{Q-learning} introduced by \citet{Watkins1989-Q-learning}. It is based on a quality function, or Q-function, defined for a control policy $v$ as follows:

\begin{equation} \label{eq:Q-fnc}
	Q_{v}(x_k, v(x_k)) = r(x_k, v(x_k) ) + J_{v}(x_{k+1}), \forall x_k \in X.	
\end{equation}

Let $\bar{Q}(x_k, v(x_k))$ denote the following function

\begin{equation} \label{eq:Q-fnc-semiopt}
	r(x_k, v(x_k) ) + J^*(x_{k+1}), \forall x_k \in X.	
\end{equation}

Then, the Bellman optimality principle (\ref{eq:DP-principle}) can be rewritten as

\begin{equation} \label{eq:Q-fnc-DP}
	J^*(x_k) = \min_v \bar{Q}( x_k, v(x_k) ), \forall x_k \in X.	
\end{equation}

with the optimal control policy (cf. (\ref{eq:opt-ctrl-policy}))

\begin{equation} \label{eq:opt-ctrl-policy-Q}
	v^*(x_k) = \arg \min_v \bar{Q}(x_k, v(x_k) ), \forall x_k \in X.	
\end{equation}

Due to the fact that Q-learning does not require the system gradients with respect to the control, it was chosen as the basic ADP control scheme in the current work. The concrete implementation of Q-learning used is provided in Sec. \ref{sec:Q-learning}. Even though the system gradients with respect to the control are not required, a means of predicting the next state of the system is still needed. Since neural networks might be computationally expensive in some applications, a simple scheme of predicting the system state based on Kalman filter \citep{Kalman1961} is suggested. The details are given in Sec. \ref{sec:KF}. A Kalman filter also helps address measurement uncertainty which is present in applications. Furthermore, it is proposed to stack the Q-functions (\ref{eq:Q-fnc-semiopt}) over a finite number of time steps to improve the policy update step. The details will be given in Sec. \ref{sec:sADP-Q}. In contrast to optimal control schemes with finite-horizon cost functions, where the running costs $r$ are stacked over a finite horizon -- such as model predictive control (MPC) -- the Q-functions are stacked in this case. Since the Q-function represents infinite-horizon optimal control, the described control scheme may be roughly seen as a mixture of infinite- and finite-horizon optimal control. For an extensive description of MPC, please refer, for example, to \citet{Garcia1989-MPC-survey}, \citet{Camacho1999-MPC} or \citet{Borrelli2011-MPC}. An insightful survey on ADP, MPC and their interrelation may be found in \citep{Bertsekas2005-ADP-MPC-survey}. 

To summarize, the contributions of this work are the following: the policy update of ADP is performed using stacked Q-function approximants; the parameters of the system model are estimated online using a Kalman filter. It is demonstrated the performance improvement is possible via the new stacked approach. The details of it are given in Sec. \ref{sec:sADP-Q} followed by comparison and case studies.

\section{Q-function algorithm implementation} \label{sec:Q-learning}

The equations (\ref{eq:Q-fnc-DP}) and (\ref{eq:opt-ctrl-policy-Q}) describe the ideal Q-learning algorithm in which the optimizing control policy must be found exactly. In the general case of nonlinear system dynamics (\ref{eq:sys-general}), it is not possible and approximation methods are thus required. In the following, all vectors are treated as column vectors if not stated otherwise. As an approximator for the Q-function, a second-order parametric form

\begin{equation} \label{eq:Q-fnc-approx-model}
	\hat{Q}(x, u) = W^T \varphi(x, u). 
\end{equation}

was chosen, where $W$ denotes the parameter vector and $\varphi$ is the regressor defined for any vectors $x, u$ as follows:

\begin{equation} \label{eq:regressor}
	\varphi(x, u) =  \mbox{vec} ( \Delta_u ( (x \vdashline u) \otimes (x \vdashline u) ) ), 
\end{equation}

where $\Delta_u ( \bullet )$ denotes the upper triangle submatrix, $\mbox{vec} ( \bullet )$ denotes the operator transforming a matrix into a vector, $(\bullet \vdashline \bullet)$ denotes stacking of vectors or matrices, and $\otimes$ denotes the Kronecker product. The approximation model (\ref{eq:Q-fnc-approx-model}) is used in linear quadratic regulators \citep{Lewis1995-opt-ctrl} where the optimal value function is known to have quadratic form. Even though, the model (\ref{eq:sys-general}) is in general nonlinear, ADP is essentially concerned with the local behavior of the system in contrast to dynamic programming. This is the reason to choose a simple approximant in the current work to ease computations. As will be shown in Sec. \ref{sec:case-study}, second-order approximant is still sufficient to provide optimizing control actions. The goal of the Q-function approximator is to adjust the parameters $W$ so as to approximately satisfy the Bellman optimality principle (\ref{eq:Q-fnc-DP}). To this end, the \emph{Bellman error} is computed:

\begin{equation} \label{eq:temp-diff}
	\begin{array}{lll}
			e_k: = & (W^+)^T \varphi( x_k, u_k^- ) - \\
			& (W^-)^T \varphi( x_{k+1}, u_k^- ) - r(x_k, u_k^-), &
	\end{array}
\end{equation}

where $W^-, u_k^-$ denote the approximant parameters and policy before policy update respectively. The parameters $W^+$ denote the new yet-to-be-found approximant parameters. They are updated iteratively. In the current work, gradient descent is applied to the optimization problem

\[
	\min \limits_{W^+} \frac{1}{2} e_k^2.
\]

This amounts to the following iterations

\begin{equation} \label{eq:W-update}
	W^+: = W^+ - \alpha e_k \varphi( x_k, u_k^- ),
\end{equation}

where $\alpha$ denotes the gradient descent gain for the \emph{critic}. The policy update is also performed by gradient descent minimizing the Q-function approximant 

\begin{equation} \label{eq:Q-fnc-approx}
	\hat{Q} := (W^+)^T \varphi( x_{k+1}, u_k^+ )
\end{equation}

with respect to $u_k^+$. The update step is thus computed as

\begin{equation} \label{eq:policy-update}
	u_k^+ := u_k^- - \beta \left( \frac{\partial \varphi( x_{k+1}, u_k ) }{\partial u_k}\Big|_{u_k = u_k^-} \right)^T  W^+,
\end{equation}

where $\alpha$ denotes the gradient descent gain for the \emph{actor}. Notice that the gradient $\frac{\partial \varphi(x, u)}{\partial u}$ for the regressor in the form (\ref{eq:regressor}) can be computed analytically. Even though the system gradient with respect to the control is not required in the described implementation of ADP, prediction of the future state $x_{k+1}$ is used in (\ref{eq:temp-diff}). In the absence of knowledge of the system dynamics (\ref{eq:sys-general}), a means of predicting $x_{k+1}$ is needed to implement Q-learning as per the method described above. To this end, a Kalman filter is used in the current work. Its details are described in the next section.

\section{Kalman--filter estimation of system gradients} \label{sec:KF}

In the current work, it is assumed that the function $f$ in (\ref{eq:sys-general}) is unknown whereas the full state vector $x$ is measured at each time step. It is suggested to consider a parametric model of $f$ and then estimate these parameters online using the knowledge of the current state $x$. A variety of methods exists for such a purpose while Kalman filter is a specific one that was designed to cope with the measurement uncertainty. In all the methods, it is crucial to use an appropriate model of the system dynamics. To this end, black-box and white-box approaches exist. The white-box approach uses a specific model structure that is derived from physical laws. In contrast, the black-box approach is based on abstract model structures. Since predicting the system state is concerned, it is important to use a dynamical model that has an ability of memorizing the system dynamics history. Recurrent neural networks may be a particular solution and are used in a number of ADP control schemes. In the current work, it is suggested to use a linear model to simplify computation. The model is formulated as follows:

\begin{equation} \label{eq:lin-pred}
	x_{k+1} = A x_k + B u_k,
\end{equation}

where $A$ and $B$ are the model parameters to be estimated. Such a model was also used in some adaptive MPC control schemes \citep{Ghaffari2013-AMPC}. The linear model (\ref{eq:lin-pred}) is simple, but has certain limitations. However, since only local behavior is essential in the current study, a linear model may be sufficient as will be demonstrated in Sec. \ref{sec:case-study}. Since the number of matrix entries in $A$ and $B$ is larger than the state vector length, it is suggested to stack the state estimates. Throughout the stack, the parameters $A$ and $B$ are assumed to stay constant. To fit the resulting estimation scheme into a Kalman filter, the following structure is used:

\[
	\begin{array}{rl}
		A_{k|k-1} & = A_{k-1|k-1} \\
		B_{k|k-1} & = B_{k-1|k-1} \\
		\left( \begin{array}{c}
			x_{k-L-1} \\
			\vdots \\
			x_k-1	
		\end{array} \right)	& = \mathcal{X}_{k-1} \mbox{vec} (A_{k-1|k-1}) + \mathcal{U}_{k-1} \mbox{vec} (B_{k-1|k-1}).
	\end{array}	
\] 

Here, $\mathcal{X}$ and $\mathcal{U}$ are matrices that formed:

\[
	\begin{array}{rl}
		\mathcal{X}_{k-1} = \left( \begin{array}{c}
			X_{k-L-1} \\
			\vdots \\
			X_{k-1}	
		\end{array} \right), & \mathcal{U}_{k-1} = \left( \begin{array}{c}
		U_{k-L-1} \\
		\vdots \\
		U_{k-1}	
		\end{array} \right)
	\end{array}
\]

with $X_i = I_n \otimes x_k^T$ and $U_i = I_n \otimes u_k^T$ accordingly. The variable $L$ denotes the stack size. Prediction of the estimate covariance is performed:

\begin{equation}
	\mathbf{P}_{k|k-1} := \mathbf{P}_{k-1|k-1} + \mathbf{Q}_k, 
\end{equation} 

where $\mathbf{Q}_k$ is the assumed state noise covariance (not to be confused with the Q-function). In the correction stage, the innovation vector is computed:

\begin{equation} 
	\varepsilon_k := \mathcal{X}_k - \mathcal{X}_{k-1} \mbox{vec} (A_{k|k-1}) + \mathcal{U}_{k-1} \mbox{vec} (B_{k|k-1}).  		
\end{equation}

Its covariance reads as:

\begin{equation} 
	\mathbf{S}_k := (\mathcal{X}_{k-1} \vdashline \mathcal{U}_{k-1}) \mathbf{P}_{k|k-1} (\mathcal{X}_{k-1} \vdashline \mathcal{U}_{k-1})^T + \mathbf{R}_k,	 		
\end{equation}  

where $\mathbf{R}_k$ is the assumed measurement noise covariance. The Kalman gain is then defined by:

\begin{equation} 
	\mathbf{K}_k := \mathbf{P}_{k|k-1} (\mathcal{X}_{k-1} \vdashline \mathcal{U}_{k-1})^T \mathbf{S}_k^{-1}.	 		
\end{equation}

The corrected estimate is computed as follows:

\begin{equation} 
	\mbox{vec} (A_{k|k} \vdashline B_{k|k})  := \mbox{vec} (A_{k|k-1} \vdashline B_{k|k-1}) + \mathbf{K}_k \varepsilon_k.	 		
\end{equation}

Finally, the corrected estimate covariance is defined by

\begin{equation} 
	\mathbf{P}_{k|k} := ( I - \mathbf{K}_k (\mathcal{X}_{k-1} \vdashline \mathcal{U}_{k-1}) ) \mathbf{P}_{k|k-1}.	 		
\end{equation}

Using the estimated parameters $A, B$, the future state can be predicted by (\ref{eq:lin-pred}) and then used in (\ref{eq:temp-diff}) of the Q-learning algorithm of Sec. \ref{sec:Q-learning}. In the next section, usage of this estimation scheme in Q-learning is discussed.

\section{Suggested algorithm} \label{sec:sADP-Q}

The Kalman filter estimation together with the algorithm of Sec. \ref{sec:Q-learning} constitute the basic ADP approach that will be denoted by ADP-Q in the following. It can be observed that the estimation scheme of Sec. \ref{sec:KF} allows predicting the state also over a stack of time steps. The idea of improving the ADP-Q is the following: computing state estimates over a stack of some $N$ time steps allows predicting the behavior of the Q-function approximant. Even though Q-function represents the optimal cost-to-go, imperfections of function approximation motivate to evaluate the control policy beyond one time step. The Q-learning algorithm of Sec. \ref{sec:Q-learning} computes the approximator parameters $W$. Using these, a stack of Q-function approximants $\hat{Q}_i = W^T \varphi(x_i, u_{i-1}), i=k+1,\dots,k+N$ is considered. The following optimization problem is addressed for the policy update:

\begin{equation} \label{eq:Q-fnc-approx-stacked} 
	\min \limits_{u_k, \dots, u_{k+N-1}} \sum_{i=k+1}^{k+N} W^T \varphi(x_i, u_{i-1}). 
\end{equation}

Denoting $\bar{u}_{k} = ( u_{k}^T \vdashline \dots \vdashline u_{k+N-1}^T )^T$, the update step (\ref{eq:policy-update}) can be thus modified as follows:

\begin{equation} \label{eq:policy-update-stacked}
	\bar{u}_{k}^+ := \bar{u}_k^- - \beta \left( \begin{array}{c}
	\sum \limits_{i=k+1}^{k+N} \left( \frac{\partial \varphi( x_i, u_{i-1} ) }{\partial u_{k}}\Big|_{\bar{u}_k^-} \right)^T  W^+ \\
	\vdots \\
	\sum \limits_{i=k+1}^{k+N} \left( \frac{\partial \varphi( x_i, u_{i-1} ) }{\partial u_{k+N-1}}\Big|_{\bar{u}_k^-} \right)^T  W^+	
	\end{array} \right),
\end{equation}

where $\bar{u}_k^+, \bar{u}_k^-$ denote the stacks of control actions of the previous and current steps respectively. Essentially, the suggested control scheme does the following: even though Q-learning is a method of infinite-horizon optimal control, due to the approximate nature of its implementation, it is worthwhile to consider the behavior of the approximant beyond one time step, i. e., over a finite stack of steps and then evaluate the control policy using (\ref{eq:Q-fnc-approx-stacked}). After updating the policy as per (\ref{eq:policy-update-stacked}), the first control action in the stack is applied which is similar to MPC. Fig. 1 briefly illustrates the described principle. According to (\ref{eq:policy-update-stacked}), the resulting control policy will thus optimize not only the current Q-function approximation, but several predicted Q-function approximations. In the following, this control scheme will be denoted by sADP-Q which stands for ``stacked'' ADP using Q-learning. Convergence of general schemes of ADP was studied by \citet{Al-Tamimi2008-VI-ADP,Heydari2014-better-prf-ADP,Liu2014-PI-ADP-conv-prf}. The recent work of \citet{Wei2016-DT-Q-learning} addressed Q-learning specifically. However, convergence of model-free ADP schemes remains a challenge as is a subject of future research.

\begin{figure} \label{fig:stacked-ADP}
	\begin{centering}
		\includegraphics{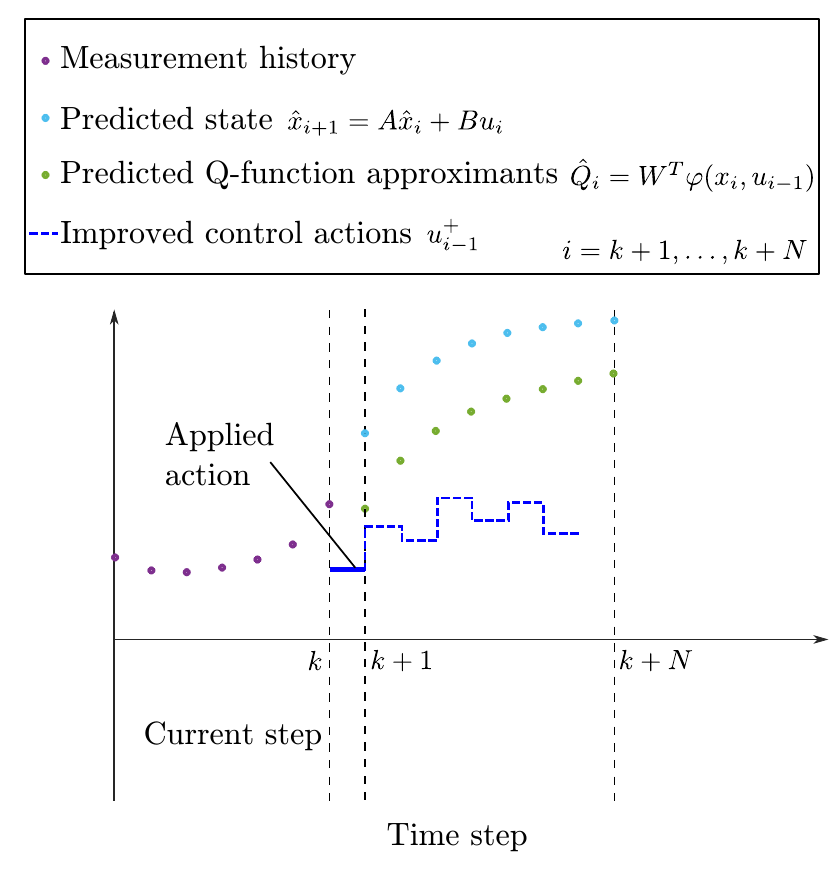}
		\par
	\end{centering}
	\caption{Prediction of Q-function approximants over a stack}	
\end{figure}

\section{Comparison with adaptive MPC} \label{sec:AMPC}

It can be noticed that a variant of an MPC algorithm without terminal cost can be obtained from the optimization problem of sADP-Q (\ref{eq:Q-fnc-approx-stacked}) by substituting $W^T \varphi(x_i, u_{i-1})$ for the running costs. The cost function gradients of the policy update (\ref{eq:policy-update-stacked}) can be computed accordingly. Table 1 briefly illustrates the similarities and key differences of sADP-Q and the said variant of MPC. It can be seen that sADP-Q introduces additional information into the optimization problem by means of the parameters $W$ along with the Q-function approximation basis functions $\varphi$.

\begin{table} \label{tab:sADP-MPC-comparison}
	\caption{Stacked Q-learning and adaptive MPC comparison}
	
	\begin{tabular}{|M{0.2\columnwidth}|M{0.3\columnwidth}|M{0.35\columnwidth}|}
		\hline		
		Control scheme & MPC & sADP-Q \tabularnewline
		\hline
		Target system & \multicolumn{2}{c|}{$x_{k+1} = f(x_k, u_k)$} \tabularnewline
		\hline
		Model identification & \multicolumn{2}{c|}{$\hat{x}_{k+1} = A \hat{x}_k + B u_k$} \tabularnewline		
		\hline
		Opt. problem & $\min \limits_{\bar{u}_k} \sum \limits_{i=k+1}^{k+N} r(x_i, u_{i-1})$ & $\min \limits_{\bar{u}_k} \sum \limits_{i=k+1}^{k+N} W^T \varphi(x_i, u_{i-1})$ \tabularnewline
		\hline
	\end{tabular}
	
\end{table}

\section{Case study} \label{sec:case-study}

The suggested control scheme was tested with the following two--dimensional nonlinear system taken from \citep[p.~68]{Lewis2013-RL-ADP}:

\begin{equation} \label{eq:exm-sys}
	x_{k+1} = \left(\begin{array}{c}
	-\sin(0.5 \cdot x_{2,k}) \\
	-\cos(1.4 \cdot x_{2,k})\sin(0.9 \cdot x_{1,k})
	\end{array}\right) + \left(\begin{array}{c}
	0 \\
	1
	\end{array}\right) u_{k}.
\end{equation}

The running cost had the following quadratic form:

\begin{equation} \label{eq:exm-run-cost}
	r(x, u) = 0.5 x^T Q x + 0.5 u^T R u
\end{equation}

with $Q = \mbox{diag}(2, 2)$ and $R = \mbox{diag}(2, 2)$.

The system was initialized at $[0.5, 1]^T$ and simulate with a step size of 1 ms. The control variable was initialized at 1. The state was measured with an additive standard normally distributed noise. The system gradients were estimated by the Kalman filter as per Sec. \ref{sec:KF}. The stack size $L$ was set to 10. The control schemes -- ADP-Q and sADP-Q -- were compared to each other and to MPC. The baseline for the comparison was represented by the gradient descent on the running cost. The baseline algorithm will be denoted by GD. It simply updates the control action by the following rule:

\begin{equation} \label{eq:GD-on-reward}
	u_k^+:= u_k^- - \beta \frac{\partial r}{\partial u}\Big|_{(x_{k+1},u_k^-)} = u_k^- - \beta (B^T x_{k+1} + u_k^-).
\end{equation}

The same parameter $\beta = 10^{-4}$ was used for all four control schemes. The ADP control schemes both used $\alpha = 0.1$ and the number of gradient descent steps equal 20 in their Q-function approximation (\ref{eq:W-update}). The approximator parameters $W$ were initalized at 1 s. The suggested control scheme -- sADP-Q -- and MPC used a stack of $N=4$ steps. First, comparison of the four control schemes is presented to show how they stabilize the state trajectory. These results can be seen in Fig. 2. It can be observed that both ADP approaches and MPC beat the GD in terms of the speed of convergence of the system to the origin whereas the sADP-Q performed better compared to its counterpart -- ADP-Q -- as well as MPC. The system dynamics model parameters $A, B$ estimated by the Kalman filter are shown in Fig. 3. As a demonstration of approximator parameters convergence, Fig. 4 illustrates $W$ during 20 steps of gradient descent adjustment by (\ref{eq:W-update}). It can be observed that after approximately 8 steps, all the parameters $W$ stabilize. The corresponding evolution of the Bellman error (\ref{eq:temp-diff}) is shown in Fig. 5. Fig. 6 shows the control actions computed by the sADP-Q. As can be seen from Fig. 7 and Fig. 8, both the running cost and the Q-function approximation were successfully minimized by the controller in about 3 s.

\begin{figure} \label{fig:state}
	\begin{centering}
		\includegraphics{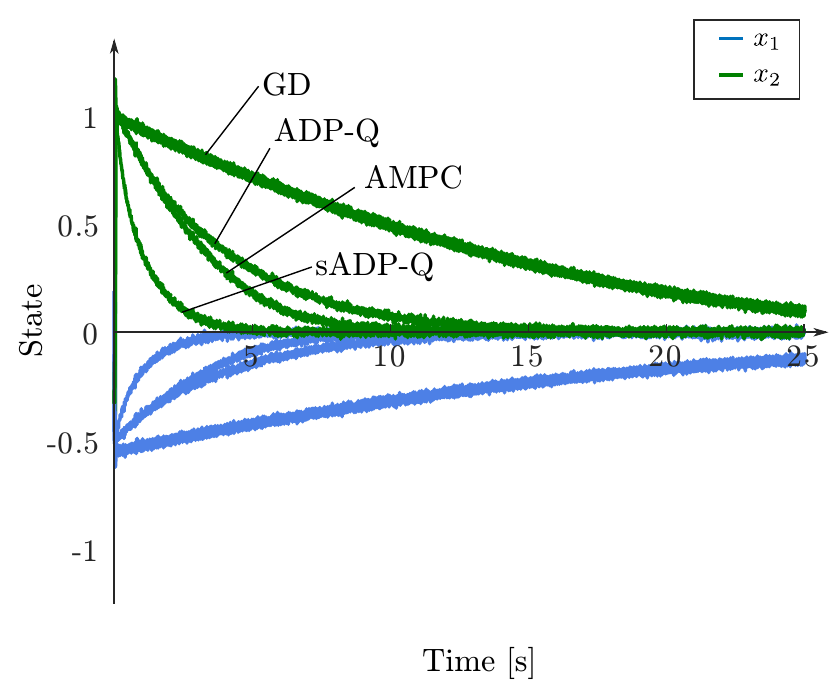}
		\par
	\end{centering}
	\caption{State trajectory under different control schemes}	
\end{figure}

\begin{figure} \label{fig:model-est-pars}
	\begin{centering}
		\includegraphics{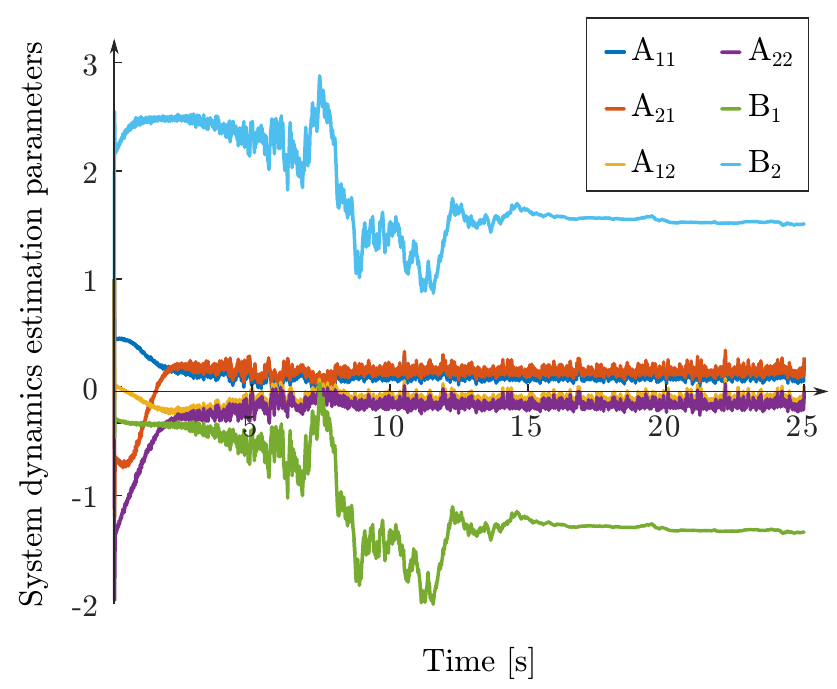}
		\par
	\end{centering}
	\caption{Estimation of the parameters of system dynamics using Kalman filter}	
\end{figure} 

\begin{figure} \label{fig:pars}
	\begin{centering}
		\includegraphics{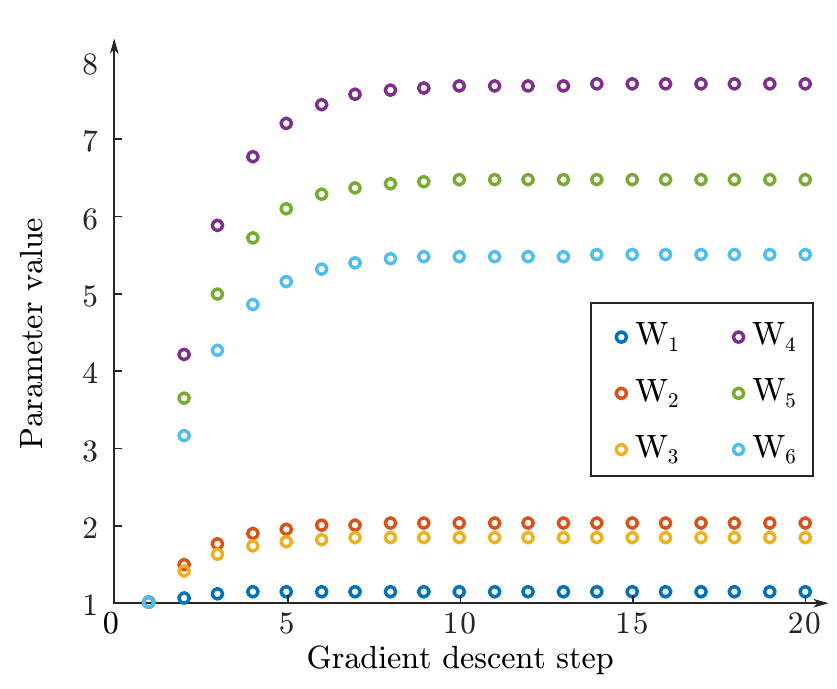}
		\par
	\end{centering}
	\caption{Convergence of the approximator parameters}	
\end{figure}

\begin{figure} \label{fig:temp-diff}
	\begin{centering}
		\includegraphics{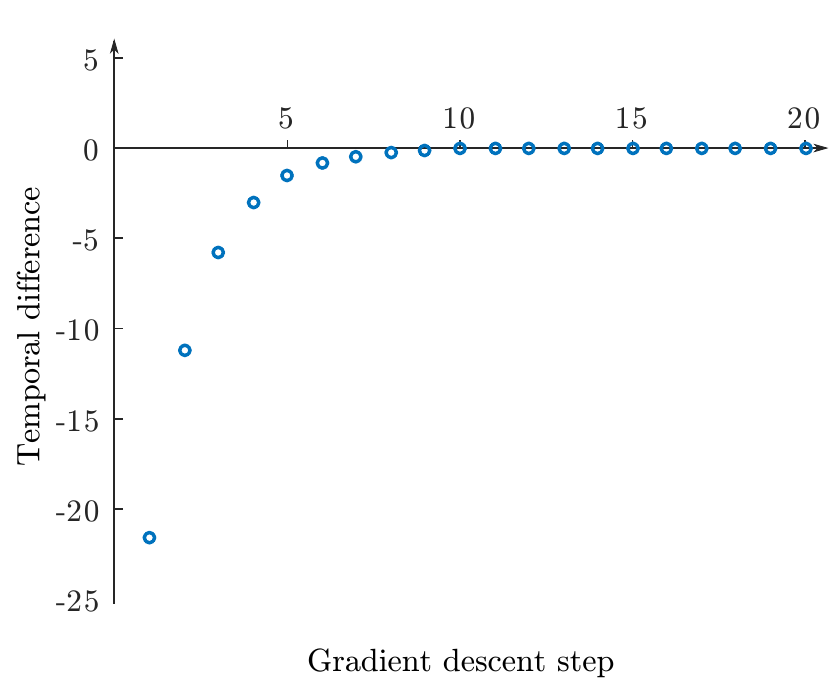}
		\par
	\end{centering}
	\caption{Convergence of the Bellman error}	
\end{figure}

\begin{figure} \label{fig:ctrl}
	\begin{centering}
		\includegraphics{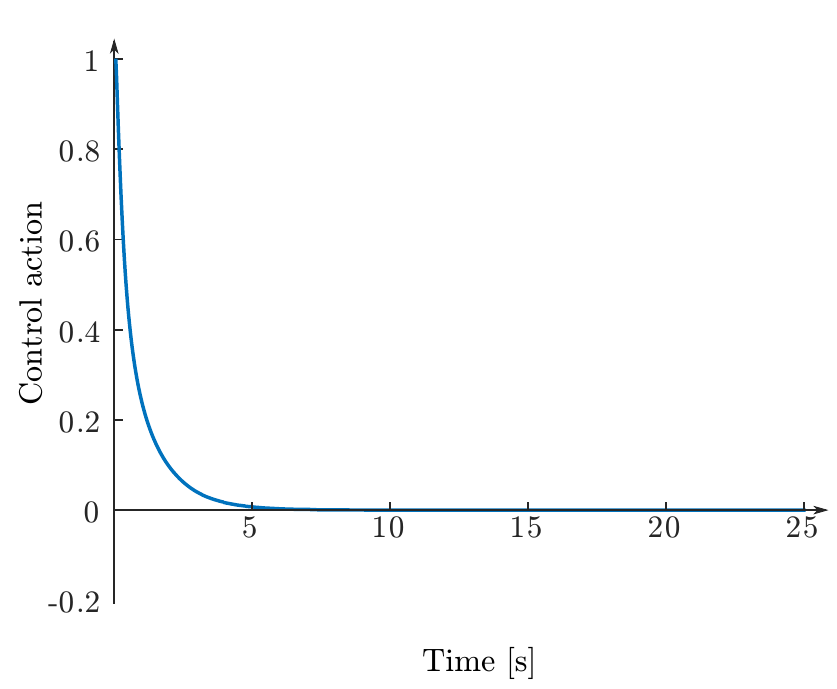}
		\par
	\end{centering}
	\caption{Control actions of sADP-Q}	
\end{figure}

\begin{figure} \label{fig:reward}
	\begin{centering}
		\includegraphics{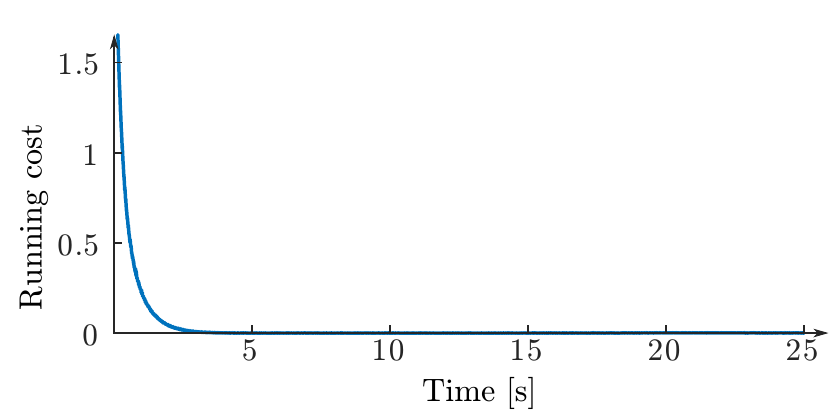}
		\par
	\end{centering}
	\caption{Running cost under sADP-Q}	
\end{figure}

\begin{figure} \label{fig:Q-fncs}
	\begin{centering}
		\includegraphics{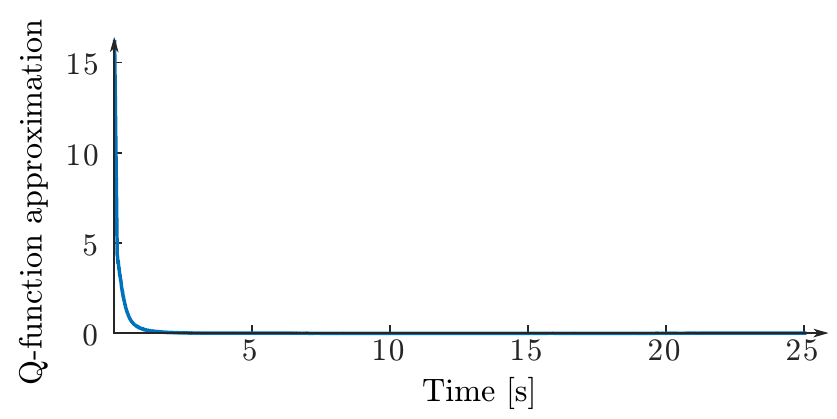}
		\par
	\end{centering}
	\caption{Q-functions approximation under sADP-Q}	
\end{figure}

\section{Discussion} \label{sec:discussion}

Since ADP is a complex algorithm that may use sophisticated function approximators, it may have a number of inappropriate settings of the tuning parameters of which may drastically deteriorate the controller performance. For instance, if a neural network is used as an approximator, it might not be clear as to how many layers and of which number of neurons to choose. In the current work, neural networks were not used, but there is still a number of tuning parameters whose proper setting might require trial-and-error. Further investigations are required as to how to optimally pick the stack length and the rest of the tuning parameters such as the gradient descent gains, or structure of the approximator. This remains an open question. The current study was dedicated to a preliminary investigation of stacked ADP and it was shown that such an approach had a certain merit. Better performance of sADP-Q compared to a variant of MPC might have resulted from the approximator model $W^T \varphi(x, u)$ in fact estimating the infinite sum of running costs which is larger than any finite sub-sum. An apparent effect of this might have been a larger magnitude of the cost function gradient. Perhaps, a similar effect in the MPC case might be achieved by an adjustment of the running cost function itself. Future research may reveal the details of stacked ADP and help understand how MPC may benefit from introducing certain features of ADP.

\section{Conclusion}

The work was concerned with an infinite-horizon optimal control of discrete systems with unknown model. As a model parameter estimation scheme, Kalman filter was used. The adaptive dynamic programming controller was based on the Q-function approach and used a simple second-order approximator. The approximator parameters were adjusted by gradient descent. To improve the controller performance, a new stacked adaptive dynamic programming was suggested. The functioning of the designed controller was compared to the basic adaptive dynamic programming, to a variant of adaptive MPC, and to gradient descent on the running cost. The results of a case study demonstrated certain benefits of the newly suggested control scheme. 

\begin{ack}
The work was supported by the Federal Ministry of Food and Agriculture (BMEL) based on a decision of the Parliament of the Federal Republic of Germany via the Federal Office for Agriculture and Food (BLE) under the innovation support programme.
\end{ack}

\bibliography{bib/DP-ADP-RL,bib/MPC,bib/opt-ctrl,bib/Kalman-filter,bib/gradient-est}                                                                                     

\begin{thebibliography}{24}
\providecommand{\natexlab}[1]{#1}
\providecommand{\url}[1]{\texttt{#1}}
\providecommand{\urlprefix}{URL }
\expandafter\ifx\csname urlstyle\endcsname\relax
  \providecommand{\doi}[1]{doi:\discretionary{}{}{}#1}\else
  \providecommand{\doi}{doi:\discretionary{}{}{}\begingroup
  \urlstyle{rm}\Url}\fi

\bibitem[{Abu-Khalaf and Lewis(2005)}]{Abu2005-ADP}
Abu-Khalaf, M. and Lewis, F.L. (2005).
\newblock Nearly optimal control laws for nonlinear systems with saturating
  actuators using a neural network {HJB} approach.
\newblock \emph{Automatica}, 41(5), 779--791.

\bibitem[{{A}l {T}amimi et~al.(2008){A}l {T}amimi, {L}ewis, and {A}bu
  {K}halaf}]{Al-Tamimi2008-VI-ADP}
{A}l {T}amimi, A., {L}ewis, F.L., and {A}bu {K}halaf, M. (2008).
\newblock {D}iscrete-{T}ime {N}onlinear {HJB} {S}olution {U}sing {A}pproximate
  {D}ynamic {P}rogramming: {C}onvergence {P}roof.
\newblock \emph{IEEE Transactions on Systems, Man, and Cybernetics, Part B
  (Cybernetics)}, 38(4), 943--949.

\bibitem[{{B}alakrishnan et~al.(2008){B}alakrishnan, {D}ing, and
  {L}ewis}]{Balakrishnan2008-ADP-survey}
{B}alakrishnan, S., {D}ing, J., and {L}ewis, F.L. (2008).
\newblock {I}ssues on stability of {ADP} feedback controllers for dynamical
  systems.
\newblock \emph{IEEE Transactions on Systems, Man, and Cybernetics, Part B:
  Cybernetics}, 38(4), 913--917.

\bibitem[{Bellman(1957)}]{Bellman1957-DP}
Bellman, R. (1957).
\newblock \emph{Dynamic Programming}.
\newblock Princeton University Press, 1st edition.

\bibitem[{Bertsekas(2005)}]{Bertsekas2005-ADP-MPC-survey}
Bertsekas, D.P. (2005).
\newblock Dynamic programming and suboptimal control: A survey from {ADP} to
  {MPC}.
\newblock \emph{European Journal of Control}, 11(4), 310--334.

\bibitem[{{B}ertsekas and {T}sitsiklis(1995)}]{Bertsekas1995-NDP}
{B}ertsekas, D.P. and {T}sitsiklis, J.N. (1995).
\newblock {N}euro-dynamic programming: an overview.
\newblock In \emph{Proceedings of the 34th IEEE Conference on Decision and
  Control}, volume~1, 560--564.

\bibitem[{Blackwell(1965)}]{Blackwell1965-DP}
Blackwell, D. (1965).
\newblock {D}iscounted dynamic programming.
\newblock \emph{The Annals of Mathematical Statistics}, 226--235.

\bibitem[{Borrelli et~al.(2011)Borrelli, Bemporad, and
  Morari}]{Borrelli2011-MPC}
Borrelli, F., Bemporad, A., and Morari, M. (2011).
\newblock Predictive {C}ontrol {F}or {L}inear and {H}ybrid {S}ystems.
\newblock \emph{Cambridge University Press}, 20, 2011.

\bibitem[{Camacho et~al.(1999)Camacho, Bordons, and
  Normey-Rico}]{Camacho1999-MPC}
Camacho, E.F., Bordons, C., and Normey-Rico, J.E. (1999).
\newblock \emph{Model Predictive Control}.
\newblock Springer.

\bibitem[{{F}errari et~al.(2011){F}errari, {S}arangapani, and
  {L}ewis}]{Ferrari2011-ADP-survey}
{F}errari, S., {S}arangapani, J., and {L}ewis, F.L. (2011).
\newblock {S}pecial issue on approximate dynamic programming and reinforcement
  learning.
\newblock \emph{Journal of Control Theory and Applications}, 9(3), 309--309.

\bibitem[{Garcia et~al.(1989)Garcia, Prett, and Morari}]{Garcia1989-MPC-survey}
Garcia, C.E., Prett, D.M., and Morari, M. (1989).
\newblock Model predictive control: theory and practice--a survey.
\newblock \emph{Automatica}, 25(3), 335--348.

\bibitem[{Ghaffari et~al.(2013)Ghaffari, Khodayari, Salehinia, Nouri-Khajavi,
  and Tafti}]{Ghaffari2013-AMPC}
Ghaffari, A., Khodayari, A., Salehinia, S., Nouri-Khajavi, M., and Tafti, M.
  (2013).
\newblock Model predictive control system design using armax identification
  method for car-following behavior.
\newblock \emph{Iranian Journal of Mechanical Engineering Transactions}, 20(1),
  48--71.

\bibitem[{{H}eydari(2014)}]{Heydari2014-better-prf-ADP}
{H}eydari, A. (2014).
\newblock {R}evisiting approximate dynamic programming and its convergence.
\newblock \emph{IEEE Transactions on Cybernetics}, 44(12), 2733--2743.

\bibitem[{Kalman and Bucy(1961)}]{Kalman1961}
Kalman, R. and Bucy, R. (1961).
\newblock {N}ew results in linear filtering and prediction theory.
\newblock \emph{Journal of Basic Engineering}, 83(3), 95--108.

\bibitem[{Lewis and Liu(2013)}]{Lewis2013-RL-ADP}
Lewis, F.L. and Liu, D. (2013).
\newblock \emph{Reinforcement {L}earning and {A}pproximate {D}ynamic
  {P}rogramming for {F}eedback {C}ontrol}, volume~17.
\newblock John Wiley \& Sons.

\bibitem[{Lewis and Syrmos(1995)}]{Lewis1995-opt-ctrl}
Lewis, F.L. and Syrmos, V.L. (1995).
\newblock \emph{{O}ptimal Control}.
\newblock John Wiley \& Sons.

\bibitem[{Lewis and Vrabie(2009)}]{Lewis2009-ADP}
Lewis, F.L. and Vrabie, D. (2009).
\newblock Reinforcement learning and adaptive dynamic programming for feedback
  control.
\newblock \emph{IEEE Circuits and Systems Magazine}, 9(3), 32--50.

\bibitem[{{L}iu and {W}ei(2014)}]{Liu2014-PI-ADP-conv-prf}
{L}iu, D. and {W}ei, Q. (2014).
\newblock {P}olicy {I}teration {A}daptive {D}ynamic {P}rogramming {A}lgorithm
  for {D}iscrete-{T}ime {N}onlinear {S}ystems.
\newblock \emph{IEEE Transactions on Neural Networks and Learning Systems},
  25(3), 621--634.

\bibitem[{{S}utton and {B}arto(1998)}]{Sutton1998-RL}
{S}utton, R.S. and {B}arto, A.G. (1998).
\newblock \emph{{R}einforcement {L}earning: {A}n {I}ntroduction}, volume~1.
\newblock MIT press Cambridge.

\bibitem[{Wald(1947)}]{Wald1947-decision-making}
Wald, A. (1947).
\newblock Foundations of a general theory of sequential decision functions.
\newblock \emph{Econometrica, Journal of the Econometric Society}, 279--313.

\bibitem[{Watkins(1989)}]{Watkins1989-Q-learning}
Watkins, C. (1989).
\newblock \emph{Learning from {D}elayed {R}ewards}.
\newblock Ph.D. thesis, University of Cambridge England.

\bibitem[{{W}ei et~al.(2016){W}ei, {L}ewis, {S}un, {Y}an, and
  {S}ong}]{Wei2016-DT-Q-learning}
{W}ei, Q., {L}ewis, F.L., {S}un, Q., {Y}an, P., and {S}ong, R. (2016).
\newblock {P}olicy {I}teration {A}daptive {D}ynamic {P}rogramming {A}lgorithm
  for {D}iscrete-{T}ime {N}onlinear {S}ystems.
\newblock \emph{IEEE Transactions on Cybernetics}, PP(99), 1--14.

\bibitem[{{W}erbos(1990)}]{Werbos1990-menu}
{W}erbos, P. (1990).
\newblock \emph{Neural {N}etworks for {C}ontrol: {A} {M}enu of {D}esigns for
  {R}einforcement {L}earning over {T}ime}.
\newblock MIT Press.

\bibitem[{{W}erbos(1992)}]{Werbos1992-ADP}
{W}erbos, P.J. (1992).
\newblock {A}pproximate dynamic programming for real-time control and neural
  modeling.
\newblock \emph{Handbook of intelligent control: Neural, fuzzy, and adaptive
  approaches}, 15, 493--525.

\end{thebibliography}
\end{document}